\documentclass[12pt]{article}
\usepackage[utf8]{inputenc}
\usepackage{amssymb,amsmath,amsfonts,eucal,mathrsfs,amsthm} 
\usepackage{hyperref}
\usepackage{empheq}
\usepackage[colorinlistoftodos]{todonotes}
\newtheorem{theorem}{Theorem}
\newtheorem{proposition}[theorem]{Proposition}
\newtheorem{lemma}[theorem]{Lemma}

\theoremstyle{definition}
\newcommand{\R}{\mathbb{R}}

\newcommand{\Sf}{\mathbb{S}}

\newcommand{\spa}{\mbox{span}}

\newcommand{\Ric}{\mbox{Ric}}

\newcommand{\rank}{\mbox{rank }}

\newcommand{\nap}{\nabla^{\perp}}
\newcommand{\nab}{\tilde\nabla}

\newcommand{\End}{\mbox{End}}

\newcommand{\trace}{\mbox{tr\,}}

\def\<{{\langle}}
\def\>{{\rangle}}

\def\a{\alpha}

\def\be{\begin{equation} }
\def\ee{\end{equation} }

\def\proof{\noindent{\it Proof:  }}
\def\qed{\ifhmode\unskip\nobreak\fi\ifmmode\ifinner
\else\hskip5 pt \fi\fi\hbox{\hskip5 pt \vrule width4 pt
height6 pt  depth1.5 pt \hskip 1pt }}
\makeatletter
\newcommand{\subjclass}[2][]{\let\@oldtitle\@title
\gdef\@title{\@oldtitle\footnotetext{#1 
\emph{Mathematics Subject Classification:} #2}}}
\newcommand{\keywords}[1]{\let\@@oldtitle\@title
\gdef\@title{\@@oldtitle\footnotetext
{\emph{Key words and phrases.} #1.}}}
\makeatother
\begin{document}

\date{}
\title{Isometric Euclidean submanifolds with\\ 
isometric Gauss maps}
\author{M. Dajczer, M. I. Jimenez and Th. Vlachos}
\keywords{Euclidean submanifold, isometric Gauss map}
\subjclass[]{53B25, 53C24, 53C42}
\maketitle

\begin{abstract}  We investigate isometric immersions 
$f\colon M^n\to\R^{n+2}$, $n\geq 3$, of Riemannian manifolds 
into Euclidean space with codimension two that admit 
isometric deformations that preserve the metric of the
Gauss map. In precise terms, the preservation of the third 
fundamental form of the submanifold must be ensured throughout 
the deformation. For minimal isometric deformations of minimal
submanifolds this is always the case. Our main result is of 
a local nature and states that if $f$ is neither minimal 
nor reducible, then it is a hypersurface of an isometrically 
deformable hypersurface $F\colon\tilde{M}^{n+1}\to\R^{n+2}$ 
such that the deformations of $F$ induce those of $f$. Moreover, 
for a particular class of such submanifolds, a complete local 
parametric description is provided.  
\end{abstract}

To what extent is a submanifold of Euclidean space locally 
determined by the properties of its Gauss map?
To lend clarity to this query, let 
$f\colon M^n\to\R^{n+p}$ be an isometric immersion of a 
Riemannian manifold of dimension $n\geq 2$ into Euclidean 
space with codimension $p$. The Gauss map of $f$ assigns 
to each point $x\in M^n$ its tangent space $f_*T_xM$ seen 
as a vector subspace of $\R^{n+p}$. Thus it is a map 
$\phi_f\colon M^n\to Gr_n(\R^{n+p})$ into the Grassmannian of 
non-oriented $n$-dimensional vector subspaces of $\R^{n+p}$.
Now suppose that there exists a non-congruent isometric immersion 
$g\colon M^n\to\R^{n+p}$ such that the Gauss maps $\phi_f$ and
$\phi_g$ of $f$ and $g$ are constrained by some predetermined 
metric condition. To address the initial question, one must 
determine the extent to which $f$ is governed by the requested
metric demand.

Certainly, the simplest requirement is to ask the Gauss maps 
$\phi_f$ and $\phi_g$ to be congruent. In simpler terms, this 
means assuming that $f$ and $g$ have the same Gauss map.
A complete answer to this question was given by Dajczer
and Gromoll in \cite{DG0}.  Roughly speaking, it was shown  
that $f$ has to be a minimal non-holomorphic isometric 
immersion of a Kaehler manifold and that $g$ must be any 
element in its one-parameter associated family of minimal 
isometric immersions.

Let $Gr_n(\R^{n+p})$ be endowed with the standard metric as a 
symmetric space.  We recall from \cite{Ob} that the pullback 
of that metric is given by the third fundamental form 
$I\!I\!I_f$ of $f$. In terms of the second fundamental form 
$\a_f\colon T_xM\times T_xM\to N_fM(x)$ of $f$ at $x\in M^n$ 
and the corresponding shape operators we have 
\be\label{tff}
I\!I\!I_f(X,Y)(x)=\sum_{i=1}^{n}\<\a_f(X,X_i),\a_f(Y,X_i)\>
=\sum_{j=1}^{p}\<A_{\xi_j}^2X,Y\>,
\ee
where $\{X_i\}_{1\leq i\leq n}$ and $\{\xi_j\}_{1\leq j\leq p}$ 
are orthonormal basis of $T_xM$ and $N_fM(x)$.

A quite weaker condition is to require the submanifolds
$f,g\colon M^n\to\R^{n+p}$ to have \emph{isometric Gauss maps}.
That means that the graphs 
$M^n\to M^n\times Gr_n(\R^{n+p})$ of $\phi_f$ and $\phi_g$
are isometric. That is, the submanifolds share identical third 
fundamental forms. This paper is devoted to examining the 
scenario that emerges under this condition. 
\vspace{1ex}

In the case of hypersurfaces, namely, for codimension $p=1$, a 
complete answer was obtained in \cite{DG1}. It was shown that the 
submanifold has to be a minimal Sbrana-Cartan hypersurface. 
We recall that Sbrana-Cartan hypersurfaces are those of rank two, 
that is, with precisely two non-zero principal curvatures 
that allow isometric deformations. For the parametric 
classification of these hypersurfaces in space forms, we refer 
to \cite{DFT} as well as to Chapter $11$ of \cite{DT} for 
additional information relevant to this paper. The case 
of surfaces with codimension $p=2$ was investigated by Vlachos 
\cite{Th}. This paper is dedicated to exploring the 
case of submanifolds of higher dimension than two but still 
in codimension $p=2$. As already evidenced by the findings 
for surfaces in \cite{Th}, preserving the metric of the Gauss 
map is quite a weak assumption. Thus, if the goal is to achieve 
parametric classifications, it becomes quite imperative to 
introduce additional assumptions. On the other hand, there are
plenty of known examples. For instance, it is quite easy to
demonstrate that any minimal submanifold allowing minimal 
isometric deformations, of which there are many, possesses 
this property. Of course, cylinders constructed over such
submanifolds are also examples. Other quite obvious examples
go as follows: Let $f_1\colon L^p\to\R^{p+1}$ be a minimal 
Sbrana-Cartan hypersurface. Given any hypersurface
$h\colon N^{n-p}\to\R^{n-p+1}$ let 
$f=h\times f_1\colon N^{n-p}\times L^p\to \R^{n+2}$
be the  extrinsic product of immersions.  The isometric
deformations to be considered are given by the ones of 
the Sbrana-Cartan hypersurface.
\vspace{1ex}

Before discussing a more elaborate class of examples, we 
recall that the \emph{relative nullity} vector subspace 
$\Delta_0(x)\subset T_xM$ at $x\in M^n$ of an isometric 
immersion $f\colon M^n\to\R^N$ is the kernel of 
its second fundamental form at that point, that is, 
$$
\Delta_0(x)=\{X\in T_xM\colon \a_f(X,Y)=0
\;\,\text{for all}\;\,Y\in T_xM\}.
$$
Also recall that on any open subset of $M^n$ where 
$\nu_f(x)=\dim\Delta_0(x)$ is constant, then 
$x\in M^n\to\Delta_0(x)$ 
is an integrable distribution whose leaves are totally geodesic 
submanifolds in $M^n$ that are mapped by $f$ onto open subsets of 
$\nu_f$-dimensional affine vector subspaces of $\R^N$.
\vspace{1ex}

Let $F\colon \tilde{M}^{n+1}\to\R^{n+2}$ be a minimal Sbrana-Cartan 
hypersurface and then let $G\colon\tilde{M}^{n+1}\to\R^{n+2}$ be an 
element of its associated one-parameter family of isometric minimal 
deformations. Let $j\colon M^n\to\tilde{M}^{n+1}$ be an isometric 
embedding whose unit normal vector field $\eta\in\Gamma(N_jM)$ 
at any point lies in the common relative nullity distribution of 
$F$ and $G$. 
It is not difficult to verify that the submanifolds 
$f=F\circ j\colon M^n\to\R^{n+2}$ and 
$g=G\circ j\colon M^n\to\R^{n+2}$ have isometric Gauss maps.
The following family of examples, fitting into this scenario,
are particularly relevant for this paper.
\vspace{1ex}

\noindent\hypertarget{Ex}{\underline{\bf{Examples}}}: 
Let $f_1\colon L^p\to\Sf^{p+1}\subset\R^{p+2}$ be a minimal 
hypersurface of rank two in the unit sphere. Then let 
$h\colon N^{n-p}\to\R^{n-p+1}$, $2\leq p\leq n-1$, be an 
immersion such that $h=(h_0,h_1)$ with $h_1>0$ for a given 
orthogonal splitting $\R^{n-p+1}=\R^{n-p}\oplus\R$. Let $M^n$
be the warped product manifold $M^n=N^{n-p}\times_{h_1}L^p$ and 
let the isometric immersion 
$f\colon M^n\to\R^{n-p}\oplus\R^{p+2}=\R^{n+2}$ be the
warped product of $h$ and $\iota\circ f_1$, where 
$\iota\colon\Sf^{p+1}\to\R^{p+2}$ denotes the inclusion. 
That is, we have that $f=(h_0,h_1f_1)$.
Given an isometric deformation $g_1\colon L^p\to\Sf^{p+1}$  
of $f_1$ then the submanifold  $g=(h_0,h_1g_1)$ is an 
isometric deformation of $f$ such that both submanifolds 
have isometric Gauss maps. 

If  $(\eta_0,\eta_1)$ is a unit normal vector field of 
$h$ then $\eta=(\eta_0,\eta_1f_1)$ is a normal vector 
field to $f$. Moreover, if $\xi$ is a unit normal vector 
field to the minimal immersion $f_1$, then seen in 
$\R^{n+2}$ it is also normal to $f$. 
Now consider the hypersurface 
$F\colon M^n\times(-\epsilon,\epsilon)\to\R^{n+2}$, 
for some $\epsilon>0$, defined by $F(x,t)=f(x)+t\eta(x)$.
Then $F$ is a minimal hypersurface of constant rank two 
having  $\xi(x,t)=\xi(x)$ as Gauss map. Moreover, the 
similarly defined minimal hypersurface 
$G\colon M^n\times(-\epsilon,\epsilon)\to\R^{n+2}$ 
is isometric to $F$. Thus $F$ and $G$ belong to a one-parameter 
family of isometric Sbrana-Cartan hypersurfaces.
\vspace{1ex}

Before stating the result of this paper, we recall 
some facts from the Sbrana-Cartan classification theory 
of isometrically deformable hypersurfaces.
\vspace{1ex}

A hypersurface $F\colon N^m\to\R^{m+1}$, $m\geq 3$, is called
\emph{surface-like} if it is a cylinder  over either a surface 
in $\R^3$ or the cone of a surface in $\Sf^3\subset\R^4$. 
According to the Sbrana-Cartan theory, any deformation of $F$ 
is given by an isometric deformation of the surface in  either
$\R^3$ or $\Sf^3$, depending on the case.

The Sbrana-Cartan hypersurface $F\colon N^m\to\R^{m+1}$, 
$m\geq 3$, is called of \emph{elliptic type} if it is not 
surface-like and the associated tensor $J\in\Gamma(\End(\Delta^\perp))$ 
satisfies $J^2=-I$; for details see \cite{DFT} or \cite{DT}.
According to the Sbrana-Cartan theory, these hypersurfaces can 
belong either to the discrete class or to the continuous class. 
The submanifolds in the former admit a single isometric deformation 
whereas the ones in the  latter a smooth one parameter family of 
isometric deformations.
\vspace{1ex}

In the following result, that $f\colon M^n\to\R^{n+2}$ is 
\emph{locally extrinsically irreducible} means that there 
is no open subset $U\subset M^n$ splitting as a Riemannian 
product  of manifolds $U=U_1\times U_2$ such that 
$f|_U$ splits extrinsically as $f|_U=f_1\times f_2$
where $f_i\colon U_i\to\R^{n_i}$ with $n+2=n_1+n_2$.
Note that this assumption excludes cylinders.
Moreover, that $f$ and $g$ \emph{extend isometrically} means 
that there is an isometric embedding $j\colon M^n\to N^{n+1}$ 
into a Riemannian manifold $N^{n+1}$ and two isometric immersions 
$F\colon N^{n+1}\to\R^{n+2}$ and $G\colon N^{n+1}\to\R^{n+2}$ 
such that $f=F\circ j$ and $g=G\circ j$.

\begin{theorem}\label{exten} Let $f\colon M^n\to\R^{n+2}$,
$n\geq 3$, be a locally extrinsically irreducible and nowhere 
minimal isometric immersion of a Riemannian manifold free of flat 
points. Let $g\colon M^n\to\R^{n+2}$ be an isometric immersion, 
not congruent to $f$ when restricted to any open subset of $M^n$, 
such that $f$ and $g$ have isometric Gauss maps.
Then along any connected component of an open and dense 
subset of $M^n$, we have the following: 
\vspace{1ex}\\
\noindent $(i)$ 
The submanifolds $f,g$ extend uniquely to isometric non-congruent 
Sbrana-Cartan hypersurfaces $F,G\colon\tilde{M}^{n+1}\to\R^{n+2}$ 
that are either surface-like or are of elliptic type. In the latter 
case, either $g$ is unique or is any element within the 
one-parameter family of isometric immersions determined by the 
isometric deformations of an elliptic Sbrana-Cartan hypersurface 
of continuous class.
\vspace{1ex}\\
\noindent $(ii)$ If the extension
$F\colon\tilde{M}^{n+1}\to\R^{n+2}$ is a minimal 
hypersurface, then either
\begin{itemize}
\item[(a)] $f$ is any submanifold given in 
\hyperlink{Ex}{\emph{Examples}} that satisfies the above 
assumptions. Therefore, it is the extrinsic warped product 
immersion of  $h\colon N^{n-p}\to\R^{n-p+1}$ and 
$\iota\circ f_1\colon L^p\to\R^{p+2}$ where 
$f_1\colon L^p\to\Sf^{p+1}$ 
is rank two minimal hypersurface, or
\item[(b)] $f=\iota\circ f_1$ where 
$f_1\colon M^n\to\Sf^{n+1}(r)$ is a minimal hypersurface 
of rank two,
\end{itemize}
and similarly for $g$ in both cases. Moreover, the hypersurface 
$F$ is the cylinder over the cone of $f_1\colon L^p\to\Sf^{p+1}$
or just the cone, and similar for $G$. 
\end{theorem}

Note that the above result shows that only some Sbrana-Cartan 
hypersurfaces carry a hypersurface whose Gauss map remains 
isometric throughout the deformations.

The submanifolds in part $(ii)$ can be characterized as the 
ones having parallel normalized mean curvature vector fields.
By normalized we mean dividing the mean curvature vector field 
by its norm to render it a unit vector field.

In the terminology introduced in \cite{DF1} the isometric 
deformations considered in the theorem are not genuine.
This means that $f$ and $g$ indeed possess isometric 
extensions, and these extensions dictate the deformations; 
see \cite{DT} for additional information. The conclusion 
of the theorem does not hold for minimal immersions. 
For instance, genuine  deformations are exemplified in 
\cite{DG2} and \cite{DV} whereas for non-genuine examples 
we refer to Theorem $5$ in \cite{DF2}. This explains 
why these submanifolds have been excluded from the 
theorem's statement.

\section{The proof}

Let $f,g\colon M^n\to\R^{n+p}$ be isometric immersions.  
The Gauss equation
$$
\<R(X,Y)Z,W\>=\<\a_f(X,W),\a_f(Y,Z)\>-\<\a_f(X,Z),\a_f(Y,W)\>
\;\,\mbox{for}\;\, 
X,Y,Z,W\in\mathfrak{X}(M)
$$
for $f$ together with the one for $g$ give that
\be\label{trGauss}
n\<\a_f(X,Y),H_f\>-I\!I\!I_f(X,Y)=\Ric(X,Y)
=n\<\a_g(X,Y),H_g\>-I\!I\!I_g(X,Y).
\ee 
Here $\Ric$ stands for the Ricci curvature tensor of $M^n$ 
whereas $H_f$ and $H_g$ denote the mean curvature vector 
fields of $f$ and $g$, respectively.
\vspace{1ex}

Assume that $f,g\colon M^n\to\R^{n+p}$ have isometric 
Gauss maps. Then \eqref{trGauss} gives
\be\label{sffH}
\<\a_f(X,Y),H_f\>=\<\a_g(X,Y),H_g\>\;\,\mbox{for any}\;\, 
X,Y\in\mathfrak{X}(M).
\ee
In particular, the mean curvature vector fields satisfy 
$\|H_f\|=\|H_g\|$.
\vspace{1ex}

In what follows, we restrict ourselves to isometric immersions 
$f,g\colon M^n\to\R^{n+2}$ with isometric Gauss maps. We have 
from \eqref{trGauss} that this condition holds true when the 
submanifolds are both minimal. Therefore, we proceed further 
under the assumption that they are non-minimal at any point.
\vspace{1ex}

In the sequel, let $\{\eta,\xi\}\subset\Gamma(N_fM)$ and 
$\{\bar{\eta},\bar{\xi}\}\subset\Gamma(N_gM)$ be local 
orthonormal normal frames where $\eta\in\Gamma(N_fM)$ and 
$\bar{\eta}\in\Gamma(N_gM)$ are vector fields colinear 
with $H_f$ and $H_g$, respectively. 
It follows from \eqref{sffH} that 
the shape operators for $f$ and $g$ satisfy 
\be\label{sheta}
A_{\eta}=\bar{A}_{\bar{\eta}}
\ee
and
\be\label{traces}
\trace A_\xi=\trace\bar{A}_{\bar{\xi}}=0.
\ee
Moreover, from \eqref{tff} and \eqref{sheta} it follows that
\be\label{xisqrd}
A_\xi^2=\bar{A}_{\bar{\xi}}^2.
\ee

To establish the validity of the theorem, we will require 
the lemmas presented next.

\begin{lemma}\label{cong}
Let $f\colon M^n\to\R^{n+2}$ be a nowhere minimal isometric
immersion free of flat points, and let $g\colon M^n\to\R^{n+2}$
be an isometric immersion such that $f$ and $g$ have isometric 
Gauss maps.
Assume that either $A_\xi=\bar{A}_{\bar{\xi}}=0$ or that
$0\neq A_\xi=\pm\bar{A}_{\bar{\xi}}$ at any point of $M^n$.
Then $f$ and $g$ are congruent submanifolds.
\end{lemma}

\proof If $A_\xi=\bar{A}_{\bar{\xi}}=0$, and considering
that $A_\eta=\bar{A}_{\bar{\eta}}$ with $\rank A_\eta\geq 2$
by assumption, then the result follows from Corollary $2.2$,
Proposition $2.9$ in \cite{DT} and the so called Fundamental 
theorem of submanifolds; for the latter see Theorem $1.10$ 
in \cite{DT}.

Assume that $0\neq A_\xi=\pm\bar{A}_{\bar{\xi}}$. Then the map
$\varphi\colon N_fM\to N_gM$ defined by $\varphi\eta=\bar\eta$ 
and $\varphi\xi=\pm\bar{\xi}$ is a vector bundle isometry 
that preserves the second fundamental form. It is then an 
elementary fact (cf.\ Lemma $4.16$ in \cite{DT}) that 
$\varphi$ is parallel, that is, it preserves the normal 
connections, and the result follows again from the Fundamental 
theorem of submanifolds.
\vspace{2ex}\qed

From the Gauss equation and \eqref{sheta} we have that
\be\label{determ}
A_\xi X\wedge A_\xi Y=\bar{A}_{\bar{\xi}}X\wedge\bar{A}_{\bar{\xi}}Y
\;\,\mbox{for any}\;\,X,Y\in\mathfrak{X}(M).
\ee
Equivalently, the symmetric bilinear map 
$\beta\colon T_xM\times T_xM\to\R^{1,1}$ defined by
$$
\beta(X,Y)=(\<A_\xi X,Y\>,\<\bar{A}_{\bar{\xi}}X,Y\>)
\;\,\mbox{for any}\;\,X,Y\in\mathfrak{X}(M)
$$
is a flat bilinear form. Here $\R^{1,1}$ is $\R^2$ endowed with 
the Lorentzian inner product
$$
\<\<(a,b),(c,d)\>\>=ac-bd
$$
and $\beta$ being flat means that
$$
\<\<\beta(X,W),\beta(Y,Z)\>\>=\<\<\beta(X,Z),\beta(Y,W)\>\>
\;\,\mbox{for all}\;\, X,Y,Z,W\in\mathfrak{X}(M).
$$
The bilinear map $\beta$ is said to be null if
$$
\<\<\beta(X,W),\beta(Y,Z)\>\>=0\;\,
\mbox{for all}\;\, X,Y,Z,W\in\mathfrak{X}(M).
$$
It is simple to check that $\beta$ is null if and only if 
$A_\xi=\pm\bar{A}_{\xi}$. 

\begin{lemma}\label{flat} Assume that $\beta$ at $x\in M^n$ 
is not null. Then the vector subspace
$$
\mathcal{N}(\beta)(x)=\{X\in T_xM\colon\beta(X,Y)
=0\;\,\mbox{for any}\;\,Y\in T_xM\}
$$
satisfies $\dim(\mathcal{N}(\beta)(x))\geq n-2$.
\end{lemma}
\proof Since $\beta$ is not null, we have that 
$\spa\{\beta(X,Y)\colon X,Y\in\mathfrak{X}(M)\}=\R^{1,1}$.
Then the claim follows from Lemma $4.20$ in \cite{DT}.\qed 

\begin{lemma}\label{rank}
Let $f,g\colon M^n\to \R^{n+2}$ be nowhere minimal
isometric immersions with isometric Gauss maps. Assume
that they are  not congruent when restricted to any
open subset of $M^n$.
If $M^n$ is free of flat points then either 
$A_\xi=\bar{A}_{\bar{\xi}}=0$ or 
$\Delta=\ker A_\xi=\ker\bar{A}_{\bar{\xi}}$ has dimension 
$n-2$.   
\end{lemma}

\proof By Lemma \ref{cong} the flat bilinear form $\beta$ 
is not null. Then Lemma \ref{flat} yields that 
$\dim(\mathcal{N}(\beta))\geq n-2$. Hence
$\mathcal{N}(\beta)\subset\ker A_{\xi}\cap\ker\bar{A}_{\bar{\xi}}$.
Thus \eqref{traces} and \eqref{xisqrd} give for $A_{\xi}$ 
and $\bar{A}_{\bar{\xi}}$ that either both vanish or that the 
common kernel $\Delta$ has dimension $n-2$.\qed
\vspace{2ex}

\noindent\emph{Proof of Theorem \ref{exten}:}
Let the connection forms $\psi,\bar{\psi}\in\Gamma(TM^*)$ 
for $f$ and $g$ be given by
$$
\nap_X\eta=\psi(X)\xi \;\;\mbox{and}\;\; 
\bar{\nabla}^\perp_X\bar{\eta}=\bar{\psi}(X)\bar\xi,
$$
where $\nap$ and $\bar{\nabla}^\perp$ are the respective
normal connections. The Codazzi equation 
$$
(\nabla_XA_\eta)Y-A_{\nap_X\eta}Y
=(\nabla_YA_\eta)X-A_{\nap_Y\eta} X
$$
and \eqref{sheta} yield that
\be\label{shxi}
\psi(X)A_\xi Y-\psi(Y)A_\xi X
=\bar{\psi}(X)\bar{A}_{\bar{\xi}}Y
-\bar{\psi}(Y)\bar{A}_{\bar{\xi}}X
\;\,\text{for any}\;\, X,Y\in\mathfrak{X}(M).
\ee

According to Lemma \ref{rank} we consider two
cases. First suppose that  $A_\xi=\bar{A}_{\bar{\xi}}=0$
along an open connected subset  $U\subset M^n$. 
Since $A_\eta=\bar{A}_{\bar{\eta}}$ and $\rank A_\eta\geq 2$,
following the arguments used in the proof of Lemma \ref{cong},
it follows that the 
submanifolds $f|_U$ and $g|_U$ are contained in  affine 
hyperplanes of $\R^{n+2}$, and thus are congruent in $\R^{n+2}$. 
But this has been ruled out. 
\vspace{1ex}
 
In the sequel, we assume that $f$ and $g$ are restricted
to a connected component of the open and dense subset where  
$\rank A_\xi=\rank\bar{A}_{\bar{\xi}}=2$, with
$\Delta=\ker A_\xi=\ker \bar{A}_{\xi}$, and that one of 
the following cases of consideration holds.
\vspace{1ex}

\noindent \emph{Case} $\psi=0$ and $\bar{\psi}\neq 0$. 
This case is not possible. In fact, from \eqref{shxi} we have
$$
\bar{\psi}(X)\bar{A}_{\bar{\xi}}Y
=\bar{\psi}(Y)\bar{A}_{\bar{\xi}}X
\;\,\mbox{for any}\;\,X,Y\in\mathfrak{X}(M).
$$
But then $\rank\bar{A}_{\bar{\xi}}<2$, which
is a contradiction. 
\vspace{1ex}

\noindent\emph{Case} $\psi\neq 0\neq\bar{\psi}$. 
First suppose that $\ker\psi=\ker\bar{\psi}$. 
If $X\in\Gamma((\ker\psi)^\perp)$ and $Y\in\Gamma(\Delta^\perp)$ 
are orthogonal then \eqref{shxi} yields
\be\label{equals}
\psi(X)A_\xi Y=\bar{\psi}(X)\bar{A}_{\bar{\xi}}Y.
\ee
By \eqref{determ} we have that $A_\xi$ and $\bar{A}_{\bar{\xi}}$ 
restricted to $\Delta^\perp$ have the same determinant. 
This together with \eqref{traces} and \eqref{equals} give  
that $A_\xi=\pm\bar{A}_{\bar{\xi}}$, and Lemma \ref{cong} 
yields a contradiction.

Therefore, we have that $\ker\psi\neq\ker\bar{\psi}$. Then 
\eqref{shxi} yields $\Delta=\ker\psi\cap\ker\bar{\psi}$.
Hence, from  \eqref{traces} and \eqref{determ} it follows that
\be\label{shxi2}
\bar{A}_{\bar{\xi}}=A_\xi\circ R,
\ee
where $R$ is an isometry of $TM$ acting as the identity on
$\Delta$.

For simplicity, henceforth we denote the restriction of any 
tensor to $\Delta^\perp$ in the same manner.
Replacing \eqref{shxi2} in \eqref{shxi} gives
\be\label{equation}
\psi(X)A_\xi Y-\psi(Y)A_\xi X
=\bar{\psi}(X)A_\xi RY-\bar{\psi}(Y)A_\xi RX
\;\,\text{for any}\;\, X,Y\in\mathfrak{X}(M).
\ee
Let $\{X,Y\}$ be an orthonormal frame of eigenvectors of 
$A_\xi$ in $\Gamma(\Delta^\perp)$ associated to the eigenvalues 
$\pm\lambda$, respectively. With respect to this frame set
\be\label{R}
R=\begin{bmatrix}
\cos\gamma & -\sin \gamma\\
\sin\gamma & \;\;\;\cos\gamma
\end{bmatrix}
\;\,\mbox{where}\;\,\gamma\in C^\infty(M).
\ee
Then we obtain from \eqref{equation} and \eqref{R} that
\be\label{psis}
\bar{\psi}=\psi\circ R.
\ee

Set $N_fM=L\oplus P$ where $L=\spa\{\xi\}$ and $P=\spa\{\eta\}$. 
Then let
\be\label{omega}
\Omega(x)
=\spa\{(\nab_X\delta)_{f_*TM\oplus P}
\colon X\in T_xM\;\text{and}\;\delta\in \Gamma(L)\}.
\ee 
Since the vector field $\xi$ is constant along $\Delta$, 
then $\Omega$ is a smooth vector subbundle of 
$f_*\Delta^\perp\oplus P$ of rank two. Let $\Lambda$ be
the line vector bundle defined by the orthogonal splitting
$f_*\Delta^\perp\oplus P=\Omega\oplus\Lambda$. Then let 
$0\neq X_0\in\Gamma(\Delta^\perp)$ be given by
$f_*X_0+\eta\in\Gamma(\Lambda)$, that is, it satisfies 
that
\be\label{xipsi1}
0=\<\nab_X(f_*X_0+\eta),\xi\>=\<A_\xi X,X_0\>+\psi(X)
\;\;\text{for any}\;\;X\in\Gamma(\Delta^\perp).
\ee

Using \eqref{shxi2}, \eqref{psis} and \eqref{xipsi1} we 
obtain that $g_*X_0+\bar{\eta}$ is orthogonal to 
$\nab_X\bar{\xi}$ for any $X\in\Delta^\perp$. Hence also
\be\label{xipsi2}
\<\bar{A}_{\bar{\xi}} X,X_0\>+\bar{\psi}(X)=0
\;\,\mbox{for any}\;\, X\in\Gamma(\Delta^\perp).
\ee

The maps $F,G\colon \tilde{M}^{n+1}
=M^n\times(-\epsilon,\epsilon)\to\R^{n+2}$
given by
\be\label{Ffnb0}
F(x,t)=f(x)+t(f_*X_0+\eta)\;\;\mbox{and}\;\;
G(x,t)=g(x)+t(g_*X_0+\bar{\eta})
\ee
define for small $\epsilon>0$ non-congruent isometric  
hypersurfaces of rank two. The fact that they are isometric
follows using \eqref{xipsi1} and \eqref{xipsi2}. Note
that the unit vector fields $\xi(x,t)=\xi(x)$ 
and $\bar{\xi}(x,t)=\bar{\xi}(x)$ are the Gauss maps of $F$ 
and $G$, respectively. Then $\Delta\oplus\spa\{\partial_t\}$
are the relative nullity subspaces of both immersions. 
\vspace{1ex}

The subsequent computations are done along $f$ and $g$, 
that is, for $t=0$. First we construct orthonormal frames 
of $\Omega$ and $\bar{\Omega}$ which are similarly defined.
Let $Y_1\in\Gamma(\Delta^\perp)$ be a unit norm vector 
field orthogonal to $X_0$. If $\rho\in C^\infty(M)$ 
is given by $\rho^2=1+1/\|X_0\|^2$, then the vector field
$Z_1=(1/\rho)((1-\rho^2)X_0+\eta)$
together with $Y_1$ form an orthonormal basis of $\Omega$.
Also $\{Y_1,\bar{Z}_1\}\in\bar{\Omega}$ is an orthonormal
frame, where $\bar{Z}_1$ is similarly defined in terms of
$\bar{\eta}$. If $A^F$ and $A^G$ denote the shape operators 
of $F$ and $G$, respectively, then 
\be\label{xipsi0}
A^FX=A_\xi X+\psi(X)\eta\;\;\mbox{and}\;\;
A^GX=\bar{A}_{\bar{\xi}} X+\bar{\psi}(X)\bar{\eta}
\;\,\mbox{for any}\;\,X\in\Gamma(\Delta^\perp).
\ee
Setting $X_1=X_0/\|X_0\|$, we obtain from \eqref{xipsi1}
and \eqref{xipsi2}  that
\be\label{rel1}
\<A_\xi X_1,X_1\>=-\|X_0\|^{-1}\psi(X_1)\;\;\mbox{and}
\;\;\<\bar{A}_{\bar{\xi}} X_1,X_1\>=-\|X_0\|^{-1}\bar{\psi}(X_1).
\ee
Then \eqref{traces} gives that
\be\label{rel4}
\<A_\xi Y_1,Y_1\>=\|X_0\|^{-1}\psi(X_1)\;\;\mbox{and}\;\;
\<\bar{A}_{\bar{\xi}} Y_1,Y_1\>=\|X_0\|^{-1}\bar{\psi}(X_1).
\ee
Using \eqref{xipsi1}, \eqref{xipsi2} and \eqref{xipsi0} 
it follows that 
\be\label{rel2}
\<A^F Y_1,Z_1\>=\rho\psi(Y_1)\;\;\mbox{and}\;\;
\<A^G Y_1,\bar{Z}_1\>=\rho\bar{\psi}(Y_1).
\ee
Since $X_0+\eta$ is contained in the common relative nullity 
subspaces of $F$ and $G$, that is, 
$A^F(X_0+\eta)=0=A^G(X_0+\eta)$,
we obtain from \eqref{xipsi1}, \eqref{xipsi2} and 
\eqref{xipsi0} that
\be\label{rel3}
\<A^FX,\eta\>=\psi(X)\;\;\mbox{and}
\;\;\<A^GX,\bar{\eta}\>=\bar{\psi}(X)
\;\,\mbox{for any}\;\,X\in\Gamma(\Delta^\perp).
\ee
Moreover, we have 
\begin{align*}
\<A^FZ_1,Z_1\>&=\frac{1}{\rho^2}\<A^F((1-\rho^2)X_0+\eta),
(1-\rho^2)X_0+\eta\>\\
&=-\frac{1}{\rho^2}\<\rho^2A^FX_0,(1-\rho^2)X_0+\eta\>
\end{align*}
and similarly for $A^G$.
It follows using \eqref{rel1} and \eqref{rel3} that
\be\label{rel5}
\<A^FZ_1,Z_1\>=-\frac{\rho^2}{\sqrt{\rho^2-1}}\psi(X_1) \;\;\mbox{and}\;\;
\<A^G\bar{Z}_1,\bar{Z}_1\>=-\frac{\rho^2}{\sqrt{\rho^2-1}}\bar{\psi}(X_1).
\ee

The frames $\{Y_1,Z_1\}$ and $\{Y_1,\bar{Z}_1\}$
coincide after identifying  $f_*TM\oplus\spa\{\eta\}$ and
$g_*TM\oplus\spa\{\bar{\eta}\}$ with the tangent space of
$M^n\times(\epsilon,\epsilon)$ at corresponding points.
Then we have from \eqref{rel4}, \eqref{rel2} and \eqref{rel5} 
that the shape operators have the expressions
\be\label{sffs}
A^F=\rho
\begin{bmatrix}
\phi\psi(X_1) & \psi(Y_1)\\
\psi(Y_1) &-\frac{1}{\phi}\psi(X_1)
\end{bmatrix},\;\;
A^G=\rho
\begin{bmatrix}
\phi\bar{\psi}(X_1) & \bar{\psi}(Y_1)\\
\bar{\psi}(Y_1) &-\frac{1}{\phi}\bar{\psi}(X_1)
\end{bmatrix}
\ee
where $\phi=\sqrt{\rho^2-1}/\rho$.
Since we can choose orientation of the frames such that
$R$ has the expression \eqref{R}, we obtain using 
\eqref{psis} that
\be\label{J}
A^G=A^F(\cos\gamma I+\sin\gamma J),
\ee
where $J$ with respect to $\{Y_1,Z_1\}$ is given by
$$
J=
\begin{bmatrix}
0 & -1/\phi\\
\phi & 0
\end{bmatrix}.
$$
Since $J^2=-I$, then either $F$ and $G$ are both surface-like or 
are isometric elliptic Sbrana-Cartan hypersurfaces.
\vspace{1ex}

Note that the second fundamental form of the surfaces involved 
in the surface-like case are as in \eqref{sffs} and satisfy the 
condition \eqref{J}. In particular, they have negative
extrinsic Gauss curvature.
\vspace{1ex}

\noindent \emph{Case} $\psi=\bar{\psi}=0$. 
The submanifolds have flat normal bundle and the shape operators 
associated to the vector fields 
$\xi$ and $\eta$ are Codazzi tensors, namely, they satisfy
\be\label{codxi}
(\nabla_XA_\xi)Y=(\nabla_YA_\xi)X\;\,\mbox{for any}\;\, 
X,Y\in\mathfrak{X}(M)
\ee
and
\be\label{codeta}
(\nabla_XA_\eta)Y=(\nabla_YA_\eta)X\;\,\mbox{for any}\;\, 
X,Y\in\mathfrak{X}(M).
\ee

Taking the $\Delta$-component of 
\eqref{codxi} for $X\in\Gamma(\Delta^\perp)$ and identity
$T\in\Gamma(\Delta)$ yields that $\Delta$ is a totally geodesic 
distribution. The Ricci equation 
$$
\<R^\perp(X,Y)\eta,\xi\>=\<[A_\eta,A_\xi]X,Y\>\;\,\mbox{for any}\;\, 
X,Y\in\mathfrak{X}(M)
$$
gives that $A_\eta$ and $A_\xi$ commute. Similarly, we have 
using \eqref{sheta} that also $A_\eta$ and $\bar{A}_{\bar{\xi}}$
commute. It follows from \eqref{traces} and 
\eqref{xisqrd} pointwise that $\bar{A}_{\bar{\xi}}=A_\xi\circ R$
being $R$ a rotation of angle $\theta$ on $\Delta^\perp$ whereas 
the identity on $\Delta$. Moreover, it follows from \eqref{codxi}
that $\theta$ is constant and  that $R$ is not the identity, up to 
sign, from Lemma \ref{cong}. Since $A_\eta$ commutes with $R$, then  
$$
A_\eta|_{\Delta^\perp}=bI\;\mbox{where}\;b\in C^\infty(M).
$$
Set $D=\ker(A_\eta-bI)$ and notice that $\Delta^\perp\subset D$. 
Hence $\dim D\geq 2$ everywhere. We assume further that we are 
restricted to  on an open subset of $M^n$ where $D$ possesses 
constant dimension. Equation \eqref{codeta} gives that
$$
\nabla_TA_\eta Z-A_\eta\nabla_TZ=Z(b)T+b\nabla_ZT-A_\eta\nabla_ZT
$$
for any $T\in\Gamma(D)$ and $Z\in\Gamma(D^\perp)$.
Then the $S$-component for any $S\in\Gamma(D)$ yields
$$
\<(A_\eta-bI) Z,\nabla_TS\>=-Z(b)\<T,S\>.
$$
It follows that the distribution $D$ is umbilical. Let  $\delta$
be the corresponding umbilical vector field.
It is easily  seen that the distribution $D$ is spherical, that 
is, $(\nabla_T\delta)_{D^\perp}=0$ for any $T\in\Gamma(D)$. 
For instance, see Exercise $1.17$ in \cite{DT}.
\vspace{1ex}

Let $U\subset M^n$ be an open subset where the totally geodesic 
distribution $\Delta$ has constant dimension being $U$ the 
saturation by maximal leaves of some cross section of the foliation.
Then the quotient space of leaves is Hausdorff and hence 
a two dimensional manifold $L^2$. Since the vector field $\xi$ 
is constant along the leaves of $\Delta$, then it determines a 
surface $h\colon L^2\to\Sf^{n+1}$ in the unit sphere. 
From now on, we restrict ourselves to open subsets where this 
situation holds. Note that these subsets form an open and dense
subset of $M^n$.

Assume that the surface $h\colon L^2\to\Sf^{n+1}$
is substantial in $\Sf^{p+1}\subset\Sf^{n+1}$ for some 
$2\leq p\leq n-1$.  This means that $p+1$ is the least dimension
of the totally geodesic sphere $\Sf^{p+1}$ that contains the 
surface. We argue that the Euclidean vector subspace $\R^{p+2}$ 
that contains $\Sf^{p+1}$ is 
spanned by the derivatives of all order of $\xi$ with respect to 
the vectors fields in $\Gamma(\Delta^\perp)$. In fact, the Euclidean 
space $\R^{p+2}$ is spanned by the derivatives of all order 
of $\xi$ in the ambient space by vectors fields of $\Gamma(TM)$. 
On one hand, we have $\nab_Z\xi=0$ for $Z\in\Gamma(\Delta)$. 
On the other hand $\nab_X\xi=-f_*A_\xi X\in \Gamma(\Delta^\perp)$ 
for $X\in\Gamma(\Delta^\perp)$. Then 
$\nab_Z\nab_X\xi\in\Gamma(\Delta^\perp)$ since the  distribution
$\Delta$ is totally geodesic. Given that the ambient 
space is flat, then 
$$
\nab_Z\nab_{X_2}\nab_{X_1}\xi
=\nab_{X_2}\nab_Z\nab_{X_1}\xi+\nab_{[Z,X_2]}\nab_{X_1}\xi
$$
for any $X_1,X_2\in\Gamma(\Delta^\perp)$.
Hence the left hand side is spanned by derivatives of $\xi$ 
along $\Gamma(\Delta^\perp)$. Similarly, we reach the same
conclusion from  
$$
\nab_Z\nab_{X_k}\dots\nab_{X_1}\xi
=\nab_{X_k}\nab_Z\dots\nab_{X_1}\xi
+\nab_{[Z,X_k]}\nab_{X_{k-1}}\dots\nab_{X_1}\xi.
$$
for any $k$ and $X_i\in\Gamma(\Delta^\perp)$, $1\leq i\leq k$.
\vspace{1ex}

Assume that $h$ is substantial in 
$\Sf^{p+1}\subset\R^{p+2}\subset\R^{n+2}$ for some  
$1\leq p\leq n-1$. We have that
$\nab_X\xi=-f_*A_\xi X\in\Gamma(\Delta^\perp)$ for 
$X\in\Gamma(\Delta^\perp)$. Then the 
derivative with respect to $Y\in\Gamma(\Delta^\perp)$ gives
\be\label{deriv}
\nab_Y\nab_X\xi
=-f_*(\nabla_YA_\xi X)_D-\<A_\xi X,A_\xi Y\>\xi
-\<A_\xi X,Y\>(b\eta+f_*\delta).    
\ee

We analyze the nature of the derivatives along vectors in 
$\Delta^\perp$ of the terms on the right hand side of 
\eqref{deriv}. Since $D$ is an umbilical distribution, we have 
that $\nab_Xf_*Z\in\Gamma (D)$ for any $X\in\Gamma(\Delta^\perp)$
and $Z\in\Gamma(\Delta\cap D)$. Therefore, the covariant derivative 
of the term $f_*(\nabla_YA_\xi X)_D$ along vector fields in 
$\Gamma(\Delta^\perp)$ belongs to the vector subbundle 
$f_*D\oplus\spa\{\xi,b\eta+f_*\delta\}$.
We have from \eqref{codeta} that $b$ is constant along $D$. 
Then, using that $D$ is spherical, that $\Delta^\perp\subset D$ 
and that $b$ is constant along $\Delta^\perp$, it follows that 
consecutive derivatives of $b\eta+f_*\delta$ with respect to 
vectors fields in $\Gamma(\Delta^\perp)$ also belong to 
$f_*D\oplus\spa\{\xi,b\eta+f_*\delta\}$. Since the surface $h$ 
is substantial in $\Sf^{p+1}\subset\R^{p+2}$, we have seen 
that $\R^{p+2}$ is spanned by the derivatives of all orders 
of $\xi$ with respect to the vector fields in 
$\Gamma(\Delta^\perp)$. It follows that
\be\label{subset}
\R^{p+2}\subset f_*D(x)\oplus\spa\{\xi (x),b\eta(x)+f_*\delta(x)\}
\;\;\mbox{at any}\;\;x\in M^n.
\ee

Decompose orthogonally $\R^{n+2}=\R^{n-p}\oplus\R^{p+2}$ and any 
$V\in \R^{n+2}$ as $V=V_1+V_2$ accordingly. Then $\eta=\eta_1+\eta_2$. 
Suppose that $\eta_1=0$ on an open subset $M^n$. 
Then $N_fM(x)\subset\R^{p+2}$ and thus $\R^{n-p}\subset f_*T_xM$ for 
any point in the open subset. Then $\R^{n-p}$  conforms a tangent 
totally geodesic distribution contained in the relative nullity 
subspaces of the submanifold. But then $f$ 
is a cylinder with an $(n-p)$-dimensional Euclidean factor on
an open subset of $M^n$, which is ruled out by assumption. 

Suppose that $\eta_2=0$ on an open subset $M^n$. 
Since $\eta$ belongs to $\R^{n-p}$ this is also 
the case of $\nab_X\eta=-f_*A_\eta X=-bf_*X$ for any 
$X\in\Gamma(\Delta^\perp)$. But since 
$f_*\Delta^\perp(x)\subset\R^{p+2}$, then $b=0$ and hence 
$\delta=0$. Having that $\xi\in\R^{p+2}$ it 
now follows from \eqref{subset} that
$\R^{p+2}=f_*D_2(x)\oplus\spa\{\xi(x)\}$ where 
$f_*D_2=f_*TM\cap\R^{p+2}$ and that
$\R^{n-p}=f_*K(x)\oplus\spa\{\eta(x)\}$ where 
$f_*K=f_*TM\cap\R^{n-p}$. Then $TM=D_2\oplus K$ where 
$D_2$ and $K$ are orthogonal totally geodesic 
distributions. Since $D_2\subset D$ from \eqref{subset},
then $\a_f(T,W)=0$ for any $T\in\Gamma(D_2)$ and 
$W\in\Gamma(K)$. 
Now we have from Theorem $8.4$ in \cite{DT} that the 
immersion is an extrinsic product of a hypersurface 
$f_1\colon M_1\to\R^{n-p}$ and a minimal hypersurface 
$f_2\colon M_2\to\R^{p+2}$ of rank two, but this also 
is ruled out by assumption.
 
In view of the above, let $U\subset M^n$ be a connected 
component of the open and dense subset where 
$\eta_1\neq 0\neq\eta_2$. Setting $f_*D_2=f_*TM\cap\R^{p+2}$ 
we have that $\dim D_2\geq p$.  Since $b=0$ implies $\delta=0$, 
it then follows from \eqref{subset} that $D_2\subset D$. 
Since $\eta_2$ is orthogonal to $D_2$ then 
\be\label{Rp}
\R^{p+2}=f_*D_2(x)\oplus\spa\{\xi(x),\eta_2(x)\}
\;\;\mbox{at any}\;\;x\in U
\ee
and thus $\dim D_2=p$ where $p\geq 2$ since 
$\Delta^\perp\subset D_2$. 

We claim that the distribution $D_2$ is spherical.
On one hand, we have from \eqref{Rp} that $\nab_Sf_*T\in\R^{p+2}$ 
for any $S,T\in\Gamma(D_2)$. On the other hand, we know that
$$
\nab_Sf_*T=f_*(\nabla_ST)_D+\<A_\xi S,T\>\xi+\<S,T\>(b\eta+f_*\delta),
\;\,\mbox{for any}\;\,S,T\in\Gamma(D_2).
$$
Since $\xi(x)\in\R^{p+2}$ then
$f_*(\nabla_ST)_D+\<S,T\>(b\eta+f_*\delta)\in\R^{p+2}$.
Thus
\be\label{exp1}
(f_*(\nabla_ST)_D)_1+\<S,T\>(b\eta_1+(f_*\delta)_1)=0
\;\,\mbox{for any}\;\,S,T\in\Gamma(D_2).
\ee
\vspace{2ex}
Using that
$$
0=\<f_*(\nabla_ST)_D,\eta\>
=\< f_*(\nabla_ST)_D,\eta_1\>+\<f_*(\nabla_ST)_D,\eta_2\>
$$
we obtain from \eqref{exp1} that
\be\label{previous}
\< f_*(\nabla_ST)_D,\eta_2\>
=\<S,T\>\<b\eta_1+(f_*\delta)_1,\eta_1\>
\;\,\mbox{for any}\;\,S,T\in\Gamma(D_2).
\ee
We have from \eqref{Rp} and \eqref{previous} that
\begin{align*}
(f_*(\nabla_ST)_D)_2&=f_*(\nabla_ST)_{D_2} 
+\|\eta_2\|^{-2}\<f_*(\nabla_ST)_D,\eta_2\>\eta_2\\
&=f_*(\nabla_ST)_{D_2}+\<S,T\>a\eta_2
\end{align*}
where $a=\|\eta_2\|^{-2}\<b\eta_1+(f_*\delta)_1,\eta_1\>$.
Then using \eqref{exp1} it follows that
\be\label{exp2}
f_*(\nabla_ST)_D-f_*(\nabla_ST)_{D_2}
=\<S,T\>f_*\sigma.
\;\,\mbox{for any}\;\,S,T\in\Gamma(D_2),
\ee
where $f_*\sigma=a\eta_2-b\eta_1-(f_*\delta)_1\in\Gamma(D)$.
Being the distribution $D$ spherical with  mean curvature 
vector field $\delta$, we have from \eqref{exp2}
that the distribution $D_2$ is umbilical
with mean curvature vector field $\tilde{\delta}=\delta+\sigma$.

Since $\nab_T\eta=-f_*A_\eta T=-bf_*T\in\R^{p+2}$ for 
any $T\in\Gamma(D_2)$ then $\nab_T\eta_1=0$. Hence, 
$\eta_1$ is parallel in $\R^{n-p}$ along $D_2$ and then 
$\|\eta_2\|$ is constant along $D_2$.

We have that
$$
\nab_Tf_*\delta=f_*\nabla_T\delta+\<A_\xi\delta,T\>\xi
+\<A_\eta T,\delta\>\eta\;\,\mbox{for any}\;\,T\in D_2,
$$
where $\nabla_T\delta$ is colinear with $T$ since $D_2\subset D$
and $D$ is a spherical distribution. 
This gives that $(f_*\delta)_1$ is parallel in $\R^{n-p}$ along 
$D_2$ and then $\|f_*\delta_2\|$ is constant along $D_2$,
where we used that $A_\xi\delta=0$ since
$\delta\in\Gamma(\Delta)$ and  that $\<A_\eta T,\delta\>
=\<bT,\delta\>=0$. It follows from
$$
0=\<f_*\sigma,\eta\>=a\|\eta_2\|^2-b\|\eta_1\|^2
-\<(f_*\delta)_1,\eta_1\>,
$$
that the function $a$ is constant along $D_2$. 
Since $\nab_T\eta=-bT$ and $\eta_1$ is constant along $D_2$ 
then $\nab_Tf_*\sigma$ is a multiple of $T$. Having that also 
$\nabla_T\delta$ is a multiple of $T$ it follows that
$\nabla_T\tilde\delta$ is a multiple of $T$ for any 
$T\in\Gamma(D_2)$. Thus the distribution $D_2$ is 
spherical as claimed.
\vspace{1ex}

We have that $f_*K=f_*TM\cap\R^{n-p}$ satisfies $\dim K= n-p-1$.
Also $\nab_Tf_*W\in\R^{n-p}$ for any $T\in\Gamma(D_2)$
and $W\in\Gamma(K)$. Thus 
$\<\nabla_TS,W\>=-\<\nabla_TW,S\>=0$ for $T,S\in\Gamma(D_2)$.
Since $\dim D=p$ hence $D_2^\perp=K\oplus\spa\{\tilde{\delta}\}$
at any point of $U$.

We claim that the distribution $D_2^\perp$ is totally 
geodesic. It holds that 
$$
\<\nab_Z\eta,f_*T\>=-\<A_\eta Z,T\>=-b\<Z,T\>=0
$$
if $Z\in\Gamma(D_2^\perp)$ and $T\in\Gamma(D_2)$.
Then $\<\nab_Z\eta_2,f_*T\>=0$ for any $Z\in\Gamma(D_2^\perp)$ 
and $T\in\Gamma(D_2)$ since $\<\nab_Z\eta_1,f_*T\>=0$
because $\nab_Z\eta_1\in\R^{n-p}$. Since $\delta$ is orthogonal 
to $D_2$ and $\xi$, we have from \eqref{Rp} that 
$(f_*\delta)_2$ and $\eta_2$ are collinear. Hence 
$\<\nab_Z(f_*\delta)_2,f_*T\>=0$ for any $Z\in\Gamma(D_2^\perp)$ 
and $T\in\Gamma(D_2)$. Since $\tilde{\delta}=\delta+\sigma$
where $f_*\sigma=a\eta_2-b\eta_1-(f_*\delta)_1$ then
$f_*\tilde{\delta}=c\eta_2-b\eta_1$ for $c\in C^\infty(U)$. 
Hence we have that $\<\nab_Zf_*\tilde{\delta},f_*T\>=0$ for any 
$Z\in\Gamma(D_2^\perp)$. On the other hand, it is trivial that
$\<\nabla_YW,T\>=0$ if $W\in\Gamma(K)$, $T\in\Gamma(D_2)$ and 
$Y\in\mathfrak{X}(U)$. We conclude that $(\nabla_ZW)_{D_2}=0$ 
for any $Z,W\in\Gamma(D_2^\perp)$, and the claim has been proved.

Finally, since $\a_f(T,Z)=0$ for any 
$T\in\Gamma(D_2)$ and $Z\in\Gamma(D_2^\perp)$, then 
Theorem $10.4$ in \cite{DT} and a result due to Nölker 
\cite{No}, which is also Theorem $10.21$ in \cite{DT}, 
gives that $f|_{U}$ is part of a warped product of 
immersions as required.\vspace{1ex}

Now assume that the surface $h\colon L^2=V/\Delta\to\Sf^{n+1}$ 
is substantial for an open subset $V\subset M^n$. 
From \eqref{subset} we have
$$
\R^{n+2}=D(x)\oplus\spa\{\xi(x),b\eta(x)+\delta(x)\}
$$
at any point $x\in V$. Then $\dim D=n$, thus $\delta(x)=0$  
and therefore $A_\eta=bI$. Note that $b(x)\neq 0$ for any 
$x\in V$ since $f$ is nowhere minimal.
Being $\eta$ parallel in the 
normal connection, we have that $f(V)$ is contained in an 
umbilical hypersurface $\Sf^{n+1}(r)\subset\R^{n+2}$. 
Note that $\xi$ is a unit normal vector field to $f|_{V}$ in 
$\Sf^{n+1}(r)$. It follows from \eqref{traces} that
$f|_V$ is a rank two minimal hypersurface of 
$\Sf^{n+1}(r)$.\vspace{1ex}

Next we argue that the submanifold $g$ is of the same type 
as $f$. The argument is for when $f$ is a warped product since 
the other case is trivial. We have seen that 
$U$ is intrinsically part of a Riemannian warped product 
determined by the distributions $D_2$ and $D_2^\perp$, where 
the former is spherical and the latter is totally geodesic. 
Since $\Delta^\perp\subset D_2$ and $\Delta$ is the common 
kernel of $A_\xi$ and $\bar{A}_{\bar{\xi}}$, it then follows 
from \eqref{sheta} that the second fundamental form of $g$ also
satisfies that $\a_g(X,Y)=0$ for any $X\in\Gamma(D_2)$ and 
$Y\in\Gamma(D_2^\perp)$. Then also $g$ is an extrinsic warped 
product of immersions. 
\vspace{1ex}

For the sake of simplicity, we continue under the assumption 
that $f$ and $g$ satisfy the requested conditions on all of 
$M^n$, rather than on a connected component within an open 
and dense subset. The maps $F,G\colon \tilde{M}^{n+1}
=M^n\times(-\epsilon,\epsilon)\to\R^{n+2}$
given by
\be\label{Ffnb}
F(x,t)=f(x)+t\eta\;\;\mbox{and}\;\;G(x,t)=g(x)+t\bar{\eta}
\ee
parametrize for some $\epsilon>0$ isometric hypersurfaces.
We know that $f$ in \hyperlink{Ex}{Examples} is a warped 
product of a hypersurface $h\colon N^{n-p}\to\R^{n-p+1}$ with 
a minimal hypersurface $f_1\colon L^p\to\Sf^{p+1}$ of rank two 
and that the Gauss map of $F$ is determined by the Gauss map 
$\xi$ of $f_1$. Since the image of $\xi$ lies in $\Sf^{p+1}$ 
then $\R^{n-p}\subset F_*T_{(x,t)}\tilde{M}$ at any
$(x,t)\in\tilde{M}$. Thus $F$ is a cylinder over a minimal 
hypersurface $F_1\colon\tilde{L}^{p+1}\to\R^{p+2}$.
From \eqref{Ffnb} we have 
$$
\<F(x,t),\xi(x,t)\>=0\;\,\mbox{at any}\;\,(x,t)\in\tilde{M},
$$
that is, the support function of $F$ vanishes. Hence $F_1$ is 
the cone over $f_1$. Thus  $F$ is a cylinder over the 
cone of a rank two minimal hypersurface
$f_1\colon L^p\to\Sf^{p+1}$, $2\leq p\leq n-1$. If $f$ is 
a rank two minimal hypersurface of a sphere and since in this 
situation the support function of $F$ also vanishes, then $F$ 
is the cone over $f$. The same assertions hold for $G$. 
In both cases, the extensions $F$ and $G$ are rank two minimal 
hypersurfaces. In particular, they are either surface-like 
or elliptic Sbrana-Cartan hypersurfaces.
\vspace{1ex}

Next we show the uniqueness of the extensions by means of 
an argument that does not rely on an assumption 
regarding the normal connection forms.
\vspace{1ex}

We consider orthonormal frames $\{\eta,\xi\}\subset\Gamma(N_fM)$ 
and $\{\bar{\eta},\bar{\xi}\}\subset\Gamma(N_gM)$ such that 
$\eta$ and $\bar{\eta}$ lie in the direction of the corresponding 
mean curvature vector fields.
Thus \eqref{sheta} and \eqref{traces} hold. We claim that, 
up to signs, $\eta$ and $\bar{\eta}$ are the unique normal vector 
fields whose corresponding shape operators coincide.
In fact, let $A_{\eta_2}=\bar{A}_{\bar{\eta}_2}$
where $\eta_2=a\eta+b\xi\in\Gamma(N_fM)$ and
$\bar{\eta}_2=c\bar{\eta}+d\bar{\xi}\in\Gamma(N_gM)$ are 
unit vector fields. Then \eqref{sheta} yields
$$
(a-c)A_\eta=d\bar{A}_{\bar{\xi}}-bA_\xi.
$$
If $a-c\neq 0$ it follows from \eqref{traces} that
$\trace A_\eta=0$. But then $f$ would be minimal, which 
has been excluded. Thus $a=c$ and it follows from
Lemma \ref{cong} that $b=d=0$, proving the claim.

Let $\bar{F},\bar{G}\colon N^{n+1}\to\R^{n+2}$ be  
non-congruent isometric extensions of $f$ and $g$. 
Let $\eta_1$ be a vector field normal to $M^n$
in $N^{n+1}$ of unit norm. Then $\bar{F}_*\eta_1$ and 
$\bar{G}_*\eta_1$ are normal to $f$ and $g$, respectively, 
and the corresponding shape operators satisfy
$A_{\bar{F}_*\eta_1}=\bar{A}_{\bar{G}_*\eta_1}$.
Hence, up to sign, we have that $\bar{F}_*\eta_1=\eta$ 
and $\bar{G}_*\eta_1=\bar{\eta}$. Thus the vector
fields $\eta$ and $\bar{\eta}$ are tangent along $f$
and $g$ to $\bar{F}(N)$ and $\bar{G}(N)$, respectively.
Hence, along $f$ and $g$ the vector fields $\xi$ and 
$\bar{\xi}$ are normal to $\bar{F}(N)$ and $\bar{G}(N)$. 
Since $\rank A_\xi=\rank \bar{A}_{\bar{\xi}}=2$ and 
$\bar{F}$ and $\bar{G}$ are non-congruent, then these
hypersurfaces have rank two. 
Therefore, there is $X_0\in\Gamma(\Delta^\perp)$, possible 
zero, such that $f_*X_0+\eta$ lies in the relative nullity of 
$\bar{F}$, because otherwise, $\bar{F}$ would have rank 
larger than two. Since the leaves of relative nullity of 
an immersion are mapped to open subsets of affine subspaces, 
then the segments $t(f_*X_0+\eta)$ are contained in 
$\bar{F}(N)$ for small $t$. Having that the same holds for 
$\bar{G}$, then these extensions coincide with $F$ and $G$ 
given either by \eqref{Ffnb0} or \eqref{Ffnb}, according to
the case.
\vspace{1ex}

To conclude the proof of the theorem, we argue regarding 
the statement there about the set of deformations with isometric 
Gauss maps.
\vspace{1ex}

Let $g_1\colon M^n\to\R^{n+2}$ be an isometric deformation
of $f$, other than $g$, with an isometric Gauss map. Then $f$ 
and $g_1$ extend uniquely to Sbrana-Cartan hypersurfaces $F_1$ 
and $G_1$. 
We have seen that $X_0\in\Gamma(\Delta^\perp)$ given by 
$$
\<A_\xi X,X_0\>+\psi(X)=0
\;\;\text{for any}\;\;X\in\Gamma(\Delta^\perp)
$$
is such that $f_*X_0+\eta$ is contained in the relative nullity 
vector subspaces of $F_1$. This shows that $F_1$ coincides with 
$F$ given by \eqref{Ffnb0} or \eqref{Ffnb}, and thus that $G_1$ 
is an isometric deformation of $F$.

If $F$ and $G$ are the Sbrana-Cartan extensions of $f$ and $g$, 
respectively, then there is an isometric embedding 
$j\colon M^n\to\tilde{M}^{n+1}$ such that $f=F\circ j$ and 
$g=G\circ j$. The second fundamental forms $A^F$ and $A^G$ 
of $F$ and $G$ are related by \eqref{J}.
If $F$ lies in the continuous class, then it admits a 
one-parameter family of deformations, say $G_s$, with
shape operators
\be\label{Js}
A^{G_s}=A^F(\cos\gamma_s I+\sin\gamma_s J)\;\,\mbox{where}\;\,
\gamma_s\in C^\infty(\tilde{M}).
\ee
Let $\xi_s$ be the Gauss map of $G_s$. Calling $g_s=G_s\circ j$, 
then the shape operator of $g_s$ associated to $\xi_s$
is given by $j_*\bar{A}_{\xi_s}=(A^{G_s}|_{j_*TM})_{j_*TM}$.
We have from \eqref{traces} that
$$
\trace ((A^{F}\circ J|_{j_*TM})_{j_*TM})=0.
$$
Therefore it follows from \eqref{Js} that $\trace\bar{A}_{\xi_s}=0$.
If $\{\eta_s,\xi_s\}\subset\Gamma(N_{g_s}M)$ is an orthonormal 
frame, then we have that $\eta_s$ lies in the direction of the 
mean curvature vector field of $g_s$. Note that $\eta_s$ is 
normal to $j$ in $\tilde{M}^{n+1}$. Therefore it follows 
from \eqref{trGauss} that $f$ and  $g_s$ have the same 
third fundamental form. Since we have seen that any isometric 
deformation of $f$ with isometric Gauss map extends to an 
isometric deformation of $F$, then the immersions $g_s$ are
all the possible isometric deformations of $f$ with isometric 
Gauss maps. If $F$ is of the discrete class, that is, if it 
admits a single isometric deformation, then also $f$ admits 
a unique isometric deformation $g$ such that $f$ and $g$ have 
isometric Gauss maps.\qed

\section{Some comments}

\noindent $(1)$ There is some additional information for the 
submanifolds in part $(ii)$ of Theorem \ref{exten}. In fact, it 
was established that the submanifold is either a warped product 
of immersions with a minimal factor or just a minimal hypersurface 
of a sphere. Then the isometric deformations in the former case 
are given in terms of the deformations of the minimal factor, 
whereas in the latter case, by the deformation in the sphere of the 
submanifold itself.
\vspace{1ex}

\noindent $(2)$ To illustrate the difficulty in constructing 
an example that belongs to part $(i)$ 
of Theorem \ref{exten} but not to part $(ii)$, an example we 
currently lack, we make the following observation.

\begin{proposition}
Let $F\colon \bar{M}^{n+1}\to \R^{n+2}$, $n\geq 3$, be an
elliptic Sbrana-Cartan hypersurface and let
$G\colon\bar{M}^{n+1}\to \R^{n+2}$ be an isometric
deformation  Let $j\colon M^n\to\bar{M}^{n+1}$ be an 
isometric immersion satisfying at any point of $M^n$ that
\be\label{ass}
\trace((A^{F}|_{j_*TM})_{j_*TM})
=\trace((A^{G}|_{j_*TM})_{j_*TM})=0.
\ee
Then $f=F\circ j\colon M^n\to\R^{n+2}$ and
$g=G\circ j\colon M^n\to\R^{n+2}$ have isometric
Gauss maps.
\end{proposition}

\proof  If $\eta$ is a unit vector field  normal to $j(M)$
in $\bar{M}^{n+1}$ and $\eta_0=F_*\eta$,
$\bar{\eta}_0=G_*\eta$, then $A_{\eta_0}=\bar{A}_{\bar{\eta}_0}$.
Let $\xi$ and $\bar{\xi}$ be unit vector fields normal
to $F$ and $G$ along $j$, respectively.
From \eqref{ass} we have  $\trace A_\xi=\trace \bar{A}_{\bar\xi}=0$.
Then from \eqref{trGauss} and \eqref{ass} it follows that $f$
and $g$ have the same third fundamental form.\qed
\vspace{2ex}

\noindent $(3)$ A very similar result to Theorem \ref{exten} 
holds if, instead of Euclidean space, we consider submanifolds 
with isometric Gauss maps in the round sphere. This similarity 
is expected since their cones are examples in Euclidean space.

\vspace{2ex}
Marcos Dajczer is partially supported by the grant 
PID2021-124157NB-I00 funded by 
MCIN/AEI/10.13039/501100011033/ `ERDF A way of making Europe',
Spain, and are also supported by Comunidad Aut\'{o}noma de la Regi\'{o}n
de Murcia, Spain, within the framework of the Regional Programme
in Promotion of the Scientific and Technical Research (Action Plan 2022),
by Fundaci\'{o}n S\'{e}neca, Regional Agency of Science and Technology,
REF, 21899/PI/22. 
\vspace{1ex}

Miguel I. Jimenez is supported by FAPESP with the grants 2022/05321-9 
and 2023/06762-1.
\vspace{1ex}

Miguel I. Jimenez expresses gratitude to the Mathematics department 
of the University of Ioannina for their kind hospitality during the 
development of this work.

\noindent Marcos Dajczer\\
Departamento de Matem\'{a}ticas\\
Universidad de Murcia,\\
E-30100 Espinardo, Murcia -- Spain\\
e-mail: marcos@impa.br
\bigskip

\noindent Miguel Ibieta Jimenez\\
Instituto de Ciências Matemáticas e de Computação\\
Universidade de São Paulo\\
São Carlos\\
SP 13566-590 -- Brazil\\
e-mail: mibieta@icmc.usp.br

\bigskip

\noindent Theodoros Vlachos\\
University of Ioannina \\
Department of Mathematics\\
Ioannina -- Greece\\
e-mail: tvlachos@uoi.gr


\begin{thebibliography}{lll}

\bibitem{DF1} M. Dajczer and L. Florit,
\emph{Genuine deformations of submanifolds},
Comm. Anal. Geom. {\bf 12} (2004), 1105--1129.

\bibitem{DF2} M. Dajczer and L. Florit,
\emph{Genuine rigidity of Euclidean submanifolds 
in codimension two},
Geom. Dedicata  {\bf 106} (2004), 195--210. 

\bibitem{DFT} M. Dajczer, L. Florit and R. Tojeiro,
\emph{On deformable hypersurfaces in space forms},
Ann. Mat. Pura Appl. {\bf 174} (1998), 361--390. 

\bibitem{DG0} M. Dajczer and D. Gromoll,
\emph{Real Kaehler submanifolds and uniqueness of the Gauss map},
J. Differential Geom. 22 (1985), 13--28.

\bibitem{DG1} M. Dajczer and D. Gromoll,
\emph{Euclidean hypersurfaces with isometric Gauss maps}, 
Math. Z. {\bf 191} (1986), 201--205. 

\bibitem{DG2} M. Dajczer and D. Gromoll, 
\emph{The Weierstrass representation for complete 
real Kaehler submanifolds of codimension two},
Invent. Math. {\bf 119} (1995), 235--242. 

\bibitem{DT} M. Dajczer and R. Tojeiro, 
``Submanifold theory. Beyond an introduction".
Universitext. Springer, New York, 2019.

\bibitem{DV} M. Dajczer and Th. Vlachos, 
\emph{A class of complete minimal submanifolds and their 
associated families of genuine deformations}, 
Comm. Anal. Geom. {\bf26} (2018), 699--721.

\bibitem{No} S. Nölker, 
\emph{Isometric immersions of warped products}, 
Differ. Geom. Appl. {\bf 6} (1996), 1--30.

\bibitem{Ob} M. Obata,
\emph{The Gauss map of immersions of Riemannian manifolds 
in spaces of constant curvature},
J. Differential Geom. {\bf 2} (1968), 217--223.

\bibitem{Th} Th. Vlachos,
\emph{Isometric deformations of surfaces preserving the 
third fundamental form},
Ann. Mat. Pura Appl. {\bf 187} (2008), 137--155.
\end{thebibliography}
\end{document}